\documentclass[12pt,a4paper]{article}
\usepackage{graphicx}
\usepackage[latin1]{inputenc}
\usepackage[T1]{fontenc}
\usepackage[english]{babel}
\usepackage{amssymb,amsmath}

\newcommand{\thmheadercommand}[1]{\textbf{\scshape{}#1}}

\def\d{{\mathrm{d}}}

\def\Z{{\mathbb{Z}}}

\renewcommand{\geq}{\geqslant}
\renewcommand{\leq}{\leqslant}

\def\eps{\varepsilon}
\renewcommand{\epsilon}{\varepsilon}
\renewcommand{\phi}{\varphi}

\DeclareMathOperator{\Id}{Id}

\newcommand{\abs}[1]{\left|\mskip1mu#1\right|}
\newcommand{\norm}[1]{\left\|#1\right\|}

\newcommand{\presgroup}[2]{\left\langle\,#1 \mid  #2\,\right\rangle}

\newcounter{prop}
\newcounter{defi}
\newcounter{thm}
\newcounter{lem}

\newenvironment{dem}[1][]{\noindent{\thmheadercommand{Proof#1}}\,\,--\,\,}{$\square$\medskip}

\newenvironment{enonce}[1]{\medskip\noindent{\thmheadercommand{#1}}\,\,--\,\,\begin{slshape}}{\end{slshape}\medskip}

\newenvironment{enonce2}[1]{\medskip\noindent{\thmheadercommand{#1}}\,\,--\,\,}{\medskip}

\newenvironment{defi}[1][]{\refstepcounter{prop}
\begin{enonce}{Definition \theprop{}#1}}{\end{enonce}}
\newenvironment{prop}[1][]{\refstepcounter{prop}
\begin{enonce}{Proposition \theprop{}#1}}{\end{enonce}}
\newenvironment{thm}[1][]{\refstepcounter{prop}
\begin{enonce}{Theorem \theprop{}#1}}{\end{enonce}}
\newenvironment{lem}[1][]{\refstepcounter{prop}
\begin{enonce}{Lemma \theprop{}#1}}{\end{enonce}}
\newenvironment{cor}[1][]{\refstepcounter{prop}
\begin{enonce}{Corollary \theprop{}#1}}{\end{enonce}}

\newenvironment{rem}[1][]{\refstepcounter{prop}
\begin{enonce2}{Remark \theprop{}#1}}{\end{enonce2}}

\newenvironment{defi*}[1][]{
\begin{enonce}{Definition#1}}{\end{enonce}}
\newenvironment{prop*}[1][]{
\begin{enonce}{Proposition#1}}{\end{enonce}}
\newenvironment{thm*}[1][]{
\begin{enonce}{Theorem#1}}{\end{enonce}}
\newenvironment{lem*}[1][]{
\begin{enonce}{Lemma#1}}{\end{enonce}}
\newenvironment{cor*}[1][]{
\begin{enonce}{Corollary#1}}{\end{enonce}}
\newenvironment{ex*}[1][]{
\begin{enonce}{Example#1}}{\end{enonce}}
\newenvironment{exo*}[1][]{
\begin{enonce2}{Exercise#1}}{\end{enonce2}}
\newenvironment{rem*}[1][]{
\begin{enonce2}{Remark#1}}{\end{enonce2}}

\title{Cogrowth and spectral gap of generic groups}

\author{Yann Ollivier}

\begin{document}
\maketitle

\begin{abstract}
We prove that that for all
$\eps$, having cogrowth exponent at most $1/2+\eps$
(in base $2m-1$ with $m$ the number of generators) is a generic property
of groups in the density model of random groups. This generalizes a
theorem of Grigorchuk and Champetier. More generally
we show that the cogrowth of a random quotient of a torsion-free
hyperbolic group stays
close to that of this group.

This proves in particular
that the spectral gap of a generic group is as large as it can be.
\end{abstract}

\paragraph{Cogrowth of generic groups.} The spectral gap of an infinite
group (with respect to a given set of generators) is a quantity
controlling the speed of convergence of the simple random walk on the
group (see~\cite{K}); up to parity problems it is equal to the first
eigenvalue of the discrete Laplacian.  By a formula of Grigorchuk
(Theorem~4.1 of~\cite{Gri},
see also section~\ref{defcgr} below) this quantity can also be expressed
combinatorially by a quantity called \emph{cogrowth}: the smaller the
cogrowth, the larger the spectral gap (see also~\cite{C}). So this is an
important quantity from the combinatorial, probabilistic and
operator-algebraic point of view (see~\cite{GdlH} or~\cite{W} and the
references therein for an overview).

In~\cite{Gri} (Theorem~7.1) and~\cite{Ch93}, Grigorchuk and Champetier
show that groups defined by a presentation satisfying the small
cancellation condition, or a weaker assumption in the case of Champetier,
with long enough relators (depending on the number of relators in in the
presentation), has a cogrowth exponent arbitrarily close to $1/2$ (the
smallest possible value), hence a spectral gap almost as large as that of
the free group with same number of generators.

We get the same conclusion for \emph{generic} groups in a precise
probabilistic meaning: that of the density model of random groups
introduced in~\cite{Gro93}, which we briefly recall in
section~\ref{defdens}.  (Note that in the density model of random groups,
if $d>0$ the number of relators is exponentially large and so
Grigorchuk's and Champetier's results do not apply).  Recall
from~\cite{Gro93} that above density $d_{\text{crit}}=1/2$, random groups
are very probably trivial.

\begin{thm}
\label{cormain}
Let $0\leq d<1/2$ be a density parameter and let $G$ be a random group on
$m\geq 2$ generators at density $d$ and length $\ell$.

Then, for any $\eps>0$, the probability that the cogrowth exponent of $G$
lies in the interval $[1/2;1/2+\eps]$ tends to $1$ as $\ell\rightarrow\infty$.
\end{thm}

In particular, this provides a new large class of groups having a large
spectral gap. 

This theorem cannot be interpreted by saying that as the relators are
very long, the geometry of the group is trivial up to scale $\ell$.
Indeed, cogrowth is an asymptotic invariant and thus takes into account
the very non-trivial geometry of random groups at scale $\ell$ (see
paragraph ``locality of cogrowth'' below). This is crudely exemplified by
the collapse of the group when density is too large.

\bigskip

Our primary motivation is the study of generic properties of groups. The
study of random groups emerged from an affirmation of Gromov
in~\cite{Gro87} that
``almost every group is hyperbolic''. Since the pioneer work of
Champetier (\cite{Ch95}) and Ol'shanski\u\i (\cite{Ols}) it has been flourishing,
now having connections with lots of topics in group theory such as
property T, the Baum-Connes conjecture, small cancellation, the
isomorphism problem, the Haagerup property, planarity of Cayley graphs...

The density model of random groups (which we recall in
section~\ref{defdens}), introduced in~\cite{Gro93}, is very rich in allowing
a precise control of the number of relators to be put in the group (and
it actually allows this number to be very large). It has proven to be
very fruitful, as random groups at different densities can have different
properties (e.g.\ property T).  See~\cite{Gh} and~\cite{Oll} for a
general discussion of random groups and the density model,
and~\cite{Gro93} for an enlightening presentation of the initial
intuition behind this model.

\paragraph{Cogrowth of random quotients.} A generic group is simply a
random quotient of a free group\footnote{There is a very interesting and
intriguing parallel approach to generic groups, developed by
Champetier in~\cite{Ch00}, which consists in considering the topological space of
all group presentations with a given number of generators. See~\cite{P}
for a description of connections of this approach with other problems in
group theory.}.  More
generally, we show that, when taking a random
quotient of a torsion-free hyperbolic group, the cogrowth of the
resulting group is very close to that of the initial group. Recall
from~\cite{Oll} that a random quotient of a torsion-free hyperbolic group
is very probably trivial above some critical density $d_{\text{crit}}$,
which precisely depends on the cogrowth of the group (see
Theorem~\ref{thmrq} in section~\ref{defdens} below).

\begin{thm}
\label{thmmain}
Let $G_0$ be a non-elementary, torsion-free hyperbolic group generated by
the elements $a_1^{\pm 1}, \ldots, a_m^{\pm 1}$. Let $\eta$ be the
cogrowth exponent of $G_0$ with respect to this generating set.

Let $0\leq d<d_{\text{crit}}$ be a density parameter and let $G$ be a
random quotient (either by plain or reduced random words) of $G_0$ at
density $d$ and length $\ell$.

Then, for any $\eps>0$, the probability that the cogrowth exponent of $G$
lies in the interval $[\eta;\eta+\eps]$ tends to $1$ when
$\ell\rightarrow\infty$.
\end{thm}

Of course Theorem~\ref{cormain} is a particular case of
Theorem~\ref{thmmain}. Also, since the cogrowth and gross cogrowth
exponent can be computed from each other by the Grigorchuk formula (see
section~\ref{defcgr}), this implies that the gross cogrowth exponent does
not change either.

This answers a very natural question arising from~\cite{Oll}: indeed, it
is known that for each torsion-free hyperbolic group, the critical
density $d_{\text{crit}}$, below which random quotients are infinite and
above which they are trivial, is equal to $1$ minus the cogrowth exponent
(resp.\ $1$ minus the gross cogrowth exponent) for a quotient by random
reduced words (resp.\ random plain words).  So wondering what happens to
the cogrowth exponent after a random quotient is very natural.

Knowing that cogrowth does not change much allows in particular to
iterate the operation of taking a random quotient. These iterated
quotients are the main ingredient in the construction by Gromov
(\cite{Gro03}) of a counter-example to the Baum-Connes conjecture with
coefficients (see also~\cite{HLS}). Without the stability of cogrowth, in
order to get the crucial cogrowth control necessary to build these
iterated quotients Gromov had to use a very indirect and non-trivial way
involving property\ T (which allows uniform control of cogrowth over all
infinite quotients of a group); this could be avoided with our argument.
So besides their interest as generic properties of groups, the results
presented here could be helpful in the field.

\begin{rem}
Theorem~\ref{thmmain} only uses the two following facts: that the random quotient
axioms of~\cite{Oll} are satisfied, and that there is a local-to-global
principle for cogrowth in the random quotient. So in particular the
result holds under slightly weaker conditions than torsion-freeness of
$G_0$, as described in~\cite{Oll} (``harmless torsion'').
\end{rem}

\paragraph{Locality of cogrowth in hyperbolic groups.} As one of our
tools we use a result about locality of cogrowth in hyperbolic groups.
Cogrowth is an asymptotic invariant, and large relations in a group can
change it noticeably. But in hyperbolic groups, if the hyperbolicity
constant is known, it is only necessary to evaluate cogrowth in some ball
in the group to get a bound for cogrowth of the group (see
Proposition~\ref{cglocal}). So in this case cogrowth is accessible to
computation.

In the case of random quotients by relators of length $\ell$, this
principle shows that it is necessary to check cogrowth up to words of
length at most $A\ell$ for some constant $A$ (which depends on density
and actually tends to infinity when $d$ is close to the critical
density), so that geometry of the quotient matters up to scale $\ell$ but
not at higher scales.

This result may have independent interest.

\paragraph{About the proofs.} The proofs make heavy use of the techniques
developed in~\cite{Ch93} and~\cite{Oll}. We hope to have included precise
enough reminders. 

As often in hyperbolic group theory, the general case is very involved
but lots of ideas are already present in the case of the free group.
So in order to help understand the structure of the argument, we first
present a proof in the case of the free group (Theorem~\ref{cormain}),
and then the proof of Theorem~\ref{thmmain} for any torsion-free
hyperbolic group.

Also, the proofs for random quotients by reduced and plain random words
are very similar. They can be treated at once using the general but
heavy terminology of~\cite{Oll}. We rather chose to present the proof of
Theorem~\ref{cormain} in the case of reduced words (for which it seems to
be more natural) and of Theorem~\ref{thmmain} in the case of plain words.

\paragraph{Acknowledgments.} I would like to thank Étienne Ghys and
Pierre Pansu for very helpful discussions and many comments on the text.
Pierre Pansu especially insisted that I should go on with this question
at a time when I had no ideas about it.
Lots of the ideas presented here emerged during my stay at the École
normale supérieure de Lyon in Spring 2003, at the invitation of Damien
Gaboriau and Étienne Ghys. I am very grateful to all the team of the
math department there for their warmth at receiving me.

\section{Definitions and notations}
\label{defs}

\subsection{Cogrowth, gross cogrowth, spectral gap}
\label{defcgr}

These are variants around the same ideas. The spectral radius of the
random walk operator on a group was studied by Kesten in~\cite{K}, and
cogrowth was defined later, simultaneously by Grigorchuk (\cite{Gri}) and
Cohen (\cite{C}).  See~\cite{GdlH} for an overview of results and open
problems about these quantities and other, related ones.

So let $G$ be an infinite group generated by the elements
$a_1^{\pm},\ldots,a_m^{\pm 1}$. Let $W_\ell$ be the set of words $w$ of
length $\ell$ in the letters $a_1^{\pm},\ldots,a_m^{\pm 1}$ such that $w$
is equal to $e$ in the group $G$. Let $W'_\ell\subset W_\ell$ be the set
of \emph{reduced} words in $W_\ell$. (Note that $W'_\ell$ is empty if $G$
is freely generated by $a_1,\ldots,a_m$.) Denote the cardinal of a set by
$\abs{.}$.

\begin{defi}[ (Cogrowth exponent)]
The \emph{cogrowth exponent} of $G$ with respect to $a_1,\ldots,a_m$ is
defined as
\[
\eta=\lim_{\substack{\ell\rightarrow\infty\\\ell\text{ even}}} \frac1\ell
\log_{2m-1} \abs{W'_\ell}
\]
or $\eta=1/2$ if $G$ is freely generated by $a_1,\ldots,a_m$.

The \emph{gross cogrowth exponent} of $G$ with respect to $a_1,\ldots,a_m$ is
defined as
\[
\theta=\lim_{\substack{\ell\rightarrow\infty\\\ell\text{ even}}} \frac1\ell
\log_{2m} \abs{W_\ell}
\]
\end{defi}

So the cogrowth exponent is the logarithm in base $2m-1$ of the cogrowth
as defined by Grigochuk and Cohen. The exponents $\eta$ and $\theta$
always lie in the interval $[1/2;1]$, with equality only in case of
$\eta$ of a free group. Amenability of $G$ is equivalent to $\eta=1$ and
to $\theta=1$.

It is shown in the references mentioned above that the limit exists. We
have to take $\ell$ even in case there are no relations of odd length in
the group (in which case $W_\ell$ is empty).

The convention for the free group is justified by the following
Grigorchuk formula (\cite{Gri}, Theorem~4.1):
\[
(2m)^\theta=(2m-1)^\eta+(2m-1)^{1-\eta}
\]
which allows to compute one exponent knowing the other (also using that
these are at least $1/2$), and shows that $\eta$ and $\theta$ vary the
same way. Given that
$\theta$ is well-defined for a free group, the formula yields
$\eta(F_m)=1/2$. As this is also the convention which makes all our
statements valid without isolating the case of a free group, we strongly
plead for this being the right convention.

The cogrowth exponent is also the exponent of growth of the kernel of the
natural map from the free group $F_m$ to $G$ sending $a_i$ to $a_i$.

The probability of return to $e$ in time $t$ of the simple random walk on
$G$ (with respect to the generators $a_1^{\pm 1},\ldots,a_m^{\pm 1}$) is
of course equal to $\abs{W_t}/ (2m)^t$. So $(2m)^{\theta-1}$ is also the
spectral radius of the random walk operator on $L^2(G)$ defined by
$Mf(x)=\frac{1}{2m} \sum f(xa_i^{\pm 1})$. This is the form studied
by Kesten (\cite{K}), who denotes by $\lambda$ this spectral radius.

Since the discrete Laplacian on $G$ is equal to the operator $\Id-M$,
$1-(2m)^{\theta-1}$ is also equal to $\min(\lambda_1,2-\lambda_\infty)$
where $\lambda_1$ is the smallest and $\lambda_\infty$ the largest
eigenvalue of the Laplacian acting on $L^2(G)$. (The problems of
$\lambda_\infty$ and of parity of $\ell$ in the definition can be avoided
by considering lazy random walks.) In particular, if $\theta$ (or $\eta$)
is small then the spectral gap $\lambda_1$ is large.

The cardinals of the sets $W_\ell$ of course satisfy the superadditivity
property $\abs{W_{\ell+\ell'}}\geq \abs{W_\ell}\abs{W_{\ell'}}$. This implies
that for any $\ell$ we have an exact (instead of asymptotic) bound
$\abs{W_\ell}\leq (2m)^{\theta\ell}$. For cogrowth this is not exactly but
almost true, due to reduction problems, and we have
$\abs{W'_{\ell+\ell'+2}}\geq \abs{W'_\ell}\abs{W'_{\ell'}}$ and the exact
inequality $\abs{W'_\ell}\leq (2m-1)^{\eta \ell+2}$. We will often
implicitly use these inequalities in the sequel.

\subsection{The density model of random groups}
\label{defdens}

A random group is a quotient of a free group $F_m=\langle
a_1,\ldots,a_m\rangle$ by (the normal closure of) a randomly chosen set
$R\subset F_m$.  Typically $R$ is viewed as a set of words in the letters
$a_i^{\pm 1}$. So defining a random group is giving a law for $R$.

More generally, given a group $G_0$ generated by the elements $a_1^{\pm
1},\ldots,a_m^{\pm 1}$, and given a set $R$ of random words in these
generators we define a random quotient of $G_0$ by $G=G_0/\langle R
\rangle$.

The density model which we now define allows a precise control on the
size of $R$: the bigger the size of $R$, the smaller the random group.
For comparison, remember the number of words of length $\ell$ in
$a_1^{\pm 1},\ldots,a_m^{\pm 1}$ is $(2m)^\ell$, and the number of
reduced words is $(2m)(2m-1)^{\ell-1}\approx (2m-1)^\ell$.

In the whole text we suppose $m\geq 2$.

\begin{defi}[ (Density model of random groups or quotients)]
Let $G_0$ be a group generated by the elements $a_1^{\pm
1},\ldots,a_m^{\pm 1}$. Let $0\leq d\leq 1$ be a density parameter.

Let $R$ be a set of $(2m)^{d\ell}$ randomly chosen words of length $\ell$
(resp.\ a set of $(2m-1)^{d\ell}$ randomly chosen reduced words of length
$\ell$), uniformly and independently picked among all those words.

We call the group $G=G_0/\langle R\rangle$ a \emph{random quotient} of
$G_0$ by plain random words (resp.\ by reduced random words), at density $d$, at length $\ell$.

In case $G_0$ is the free group $F_m$ and reduced words are taken, we
simply call $G$ a \emph{random group}.
\end{defi}

In this definition, we can also replace ``words of length $\ell$ '' by
``words of length between $\ell$ and $\ell+C$ '' for any constant $C$; the
theorems presented thereafter remain valid. In~\cite{Oll}, section~4, we
describe generalizations of these models.

The interest of the density model was established by the following
theorem of Gromov, which shows a sharp phase transition between infinity
and triviality of random groups.

\begin{thm}[ (M.~Gromov, \cite{Gro93})]
Let $d<1/2$. Then with probability tending to $1$ as $\ell$ tends to
infinity, random groups at density $d$ are infinite hyperbolic.

Let $d>1/2$. Then with probability tending to $1$ as $\ell$ tends to
infinity, random groups at density $d$ are either $\{e\}$ or $\Z/2\Z$.
\end{thm}

(The occurrence of $\Z/2\Z$ is of course due to the case when $\ell$ is
even; this disappears if one takes words of length between $\ell$ and
$\ell+C$ with $C\geq 1$.)

Basically, $d\ell$ is to be interpreted as the ``dimension'' of the
random set $R$ (see the discussion in~\cite{Gro93}). As an illustration,
if $L< 2d\ell$ then very probably there will be two relators in $R$
sharing a common subword of length $L$. Indeed, the dimension of the
couples of relators in $R$ is $2d\ell$, whereas sharing a common subword
of length $L$ amounts to $L$ ``equations'', so the dimension of those
couples sharing a subword is $2d\ell-L$, which is positive if $L<2d\ell$.
This ``shows'' in particular that at density $d$, the small cancellation
condition $C'(2d)$ is satisfied.

\bigskip

Since a random quotient of a free group is hyperbolic, one can wonder if
a random quotient of a hyperbolic group is still hyperbolic. The answer
is basically yes, and the critical density in this case is linked to the
cogrowth exponent of the initial group.

\begin{thm}[ (Y.~Ollivier, \cite{Oll})]
\label{thmrq}
Let $G_0$ be a non-elementary,
torsion-free hyperbolic group, generated by the elements $a_1^{\pm
1},\ldots,a_m^{\pm 1}$, with cogrowth exponent $\eta$ and gross cogrowth
exponent $\theta$.

Let $0\leq d \leq 1$ be a density parameter, and set
$d_{\text{crit}}=1-\theta$ (resp.\ $d_{\text{crit}}=1-\eta$).

If $d<d_{\text{crit}}$, then a random quotient of $G_0$ by plain (resp.\
reduced) random words is infinite hyperbolic, with probability tending to
$1$ as $\ell$ tends to infinity.

If $d>d_{\text{crit}}$, then a random quotient of $G_0$ by plain (resp.\
reduced) random words is either $\{e\}$ or $\Z/2\Z$, with probability tending to
$1$ as $\ell$ tends to infinity.
\end{thm}

This is the context in which Theorem~\ref{thmmain} is to be understood.

\subsection{Hyperbolic groups and isoperimetry of van Kampen diagrams}
\label{defiso}

Let $G$ be a group given by the finite presentation
$\presgroup{a_1,\ldots,a_m}{R}$. Let $w$ be a word in the $a_i^{\pm
1}$'s. We denote by $\abs{w}$ the number of letters of $w$, and by
$\norm{w}$ the distance from $e$ to $w$ in the Cayley graph of the
presentation, that is, the minimal length of a word representing the same
element of $G$ as $w$.

Let $\lambda$ be the maximal length of a relation in $R$.

\def\d{\partial}
\def\A{\mathcal{A}}

We refer to~\cite{LS} for the definition and basic properties of van
Kampen diagrams. Remember that a word represents the neutral element of
$G$ if and only if it is the boundary word of some van Kampen diagram. If
$D$ is a van Kampen diagram, we denote its number of faces by $\abs{D}$
and its boundary length by $\abs{\d D}$.

It is known (\cite{Sh}) that $G$ is hyperbolic if and only if there
exists a constant $C_1>0$ such that for any (reduced) word $w$
representing the neutral element of $G$, there exists a van
Kampen diagram with boundary word $w$, and with at most $\abs{w}/C_1$ faces.
This can be reformulated as: for any word $w$ representing the neutral
element of $G$, there exists a van Kampen diagram with boundary word $w$
satisfying the isoperimetric inequality
\[
\abs{\d D}\geq C_1 \abs{D}
\]

We are going to use a \emph{homogeneous} way to write this inequality.
The above form compares the boundary length of a van Kampen diagram to
its number of faces. This amounts to comparing a length with a number,
which is not very well-suited for geometric arguments, especially when
dealing with groups having relations of very different lengths.

So let $D$ be a van Kampen diagram w.r.t.\ the presentation and define the
\emph{area} of $D$ to be
\[
\A(D)=\sum_{f\text{ face of }D} \abs{\d f}
\]
which is also the number of external edges (not couting ``filaments'') plus twice the number of
internal ones. This has, heuristically speaking, the homogeneity of a
length.

It is immediate to see that if $D$ satisfies $\abs{\d D}\geq C_1
\abs{D}$, then we have $\abs{\d D}\geq C_1\,\A(D)/\lambda$ (recall
$\lambda$ is the maximal length of a relation in the presentation).
Conversely, if $\abs{\d D}\geq C_2\,\A(D)$, then $\abs{\d D}\geq
C_2\abs{D}$. So we can express the isoperimetric inequality using $\A(D)$
instead of $\abs{D}$.

Say a diagram is \emph{minimal} if it has minimal area for a given
boundary word. So $G$ is hyperbolic if and only if there exists a
constant $C>0$ such that every minimal van Kampen diagram satisfies the
isoperimetric inequality
\[
\abs{\d D}\geq C\,\A(D)
\]

This formulation is homogeneous, that is, it compares a length to a
length. This inequality is the one that naturally arises in $C'(\alpha)$
small cancellation theory (with $C=1-6\alpha$) as well as in random
groups at density $d$ (with $C=\frac12-d$). So in these contexts the value of
$C$ is naturally linked with some parameters of the presentation.

This kind of isoperimetric inequality is also the one appearing in the
assumptions of Champetier in~\cite{Ch93}, in random quotients of
hyperbolic groups (cf.~\cite{Oll}) and in the (infinitely presented)
limit groups constructed by Gromov in~\cite{Gro03}.  So we think this is
the right way to write the isoperimetric inequality when the lengths of
the relators are very different.

\section{Locality of cogrowth in hyperbolic groups}

The goal of this section is to show that in a hyperbolic group, in order
to estimate cogrowth (which is an asymptotic invariant), it is enough to
check only words of bounded length, where the bound depends on the
quality of the isoperimetric inequality in the group.

Everything here is valid, \emph{mutatis mutandis}, for cogrowth and gross
cogrowth.

Here $G=\presgroup{a_1,\ldots,a_m}{R}$ ($m\geq 2$) is a hyperbolic group
and $W_\ell$ is the set of reduced words of length $\ell$ in the
$a_i^{\pm 1}$ equal to $e$ in $G$. Let also $\lambda$ be the maximal
length of a relation in $R$.

As explained above, hyperbolicity of $G$ amounts to the existence of some
constant $C>0$ such that 
any minimal van Kampen diagram $D$ over this
presentation satisfies the isoperimetric inequality
\[
\abs{\d D}\geq C \A(D)
\]

\bigskip

We will prove the following.

\begin{prop}\label{cglocal}
Suppose that, for some $A>1$, for any $A\lambda/4\leq \ell\leq A\lambda$ one has
\[
\abs{W_\ell}\leq (2m-1)^{\eta \ell}
\]
for some $\eta\geq 1/2$.

Then for any $\ell\geq A\lambda/4$,
\[
\abs{W_\ell}\leq (2m-1)^{\eta\ell(1+o(1)_{A\rightarrow\infty})}
\]
where the constant implied in $o(1)$ depends only on $C$.
\end{prop}

It follows from the proof that actually $1+o(1)\leq \exp
\frac{200}{C\sqrt{A}}$, so that is it enough to take $A\approx
40000/C^2$ for a good result.

\begin{dem}

First we need some simple lemmas.

The \emph{distance to boundary} of a face of a van Kampen diagram is the
minimal length of a sequence of faces adjacent by an edge, beginning with
the given face and ending with a face adjacent to the boundary (so that a
boundary face is at distance $1$ from the boundary).
%We say that a
%diagram is $E$-narrow if each face is at distance at most $E$ from the
%boundary.

Set $\alpha=1/\log(1/(1-C))\leq 1/C$, where we can suppose $C\leq 1$.

%It can be shown that the hyperbolicity constant of such a group is at
%most $\mathrm{Cst}.\lambda/\alpha$ where $\mathrm{Cst}$ is a universal
%constant.

\begin{lem}
Let $D$ be a minimal van Kampen diagram. Then $D$ can be written as a
disjoint union $D=D_1\cup D_2$ (with maybe $D_2$ not connected)
such that each face of $D_1$ is at
distance at most $\alpha \log(\A(D)/\lambda)$ from the boundary of
$D$, and $D_2$ has area at most $\lambda$.
\end{lem}

%Note that this is less than $\alpha\lambda(\log \abs{\d D}+\log 1/C)$.

\begin{dem}
Since $D$ is minimal it satisfies the isoperimetric inequality $\abs{\d
D}\geq C \A(D)$. Thus, the cumulated area of the faces of $D$ which are
adjacent to the boundary is at least $C \A(D)$, and so the cumulated area
of the faces at distance at least $2$ is at most $(1-C)\A(D)$.

Applying the same reasoning to the (maybe not connected) diagram obtained
from $D$ by removing the boundary faces, we get by induction that the
cumulated area of the faces of $D$ lying at distance at least $k$ from
the boundary is at most $(1-C)^{k-1}\A(D)$. Taking $k=1+\alpha
\log(\A(D)/\lambda)$ (rounded up to the nearest integer) provides the
desired decomposition.
\end{dem}

In the sequel we will neglect divisibility problems (such as the length
of a diagram being a multiple of $4$).

\begin{lem}
Let $D$ be a minimal van Kampen diagram. $D$ can be partitioned into two diagrams
$D'$, $D''$ by cutting it along a path of length at most
$\lambda+2\alpha\lambda\log (\A(D)/\lambda)$ such that each of $D'$ and $D''$ contains at least
one quarter of the boundary of $D$.
\end{lem}

(Here a \emph{path} in a diagram is meant to be a path in its $1$-skeleton.)

\begin{dem}
Consider the decomposition $D=D_1\cup D_2$ of the previous lemma, and
first suppose that $D_2$ is empty, so that any face of $D_1$ lies at
distance at most $\alpha\lambda\log (\A(D)/\lambda)$ from the boundary.

Let $L$ be the boundary length of $D$ and mark four points $A,B,C,D$ on
$\d D$ at distance $L/4$ of each other. As $D$ is $\alpha
\log (\A(D)/\lambda)$-narrow, there exists a path of length at most
$2\alpha\lambda\log (\A(D)/\lambda)$ joining either a point of $AB$ to a point of $CD$
or a point of $AD$ to a point of $BC$, which provides the desired
cutting.

Now if $D_2$ was not empty, first retract each connected component of
$D_2$ to a point: the reasoning above shows that there exists a path of
length at most $2\alpha\lambda\log (\A(D)/\lambda)$ joining either a
point of $AB$ to a point of $CD$
or a point of $AD$ to a point of $BC$, \emph{not counting the length in
$D_2$}. But since the sum of the lengths of the faces of $D_2$ is at most
$\lambda$, the cumulated length of the travel in $D_2$ is at most
$\lambda$, hence the lemma.
\end{dem}

The cardinal of the $W_\ell$'s (almost in the case of cogrowth, see
above) satisfy the supermultiplicativity property $\abs{W_\ell}\geq
\abs{W_{\ell-L}}\abs{W_L}$. Using narrowness of diagrams we are able to
show a converse inequality, which will enable us to control cogrowth.

\begin{cor}\label{quasimult}
We have, up to parity problems,
\begin{eqnarray*}
\abs{W_\ell} &\leq &
%\abs{W_{\lambda}}
\sum_{\ell/4\leq \ell'\leq 3\ell/4}
\abs{W_{\ell'+2\alpha\lambda\log (\ell/C\lambda)+\lambda}}
\abs{W_{\ell-\ell'+2\alpha\lambda\log(\ell/C\lambda)+\lambda}}
\\ &\leq &
\frac{\ell}{\lambda}
%\abs{W_{\lambda}}
\max_{\ell/4\leq \ell'\leq 3\ell/4}
\abs{W_{\ell'+2\alpha\lambda\log (\ell/C\lambda)+3\lambda}}
\abs{W_{\ell-\ell'+2\alpha\lambda\log(\ell/C\lambda)+3\lambda}}
\end{eqnarray*}
\end{cor}

\begin{dem}
Any word in $W_\ell$ is the boundary word of some (minimal) van Kampen diagram $D$
with boundary length $\ell$, and
so the first inequality follows from the previous lemma, together with
the inequality $\A(D)\leq \abs{\d D}/C$.

The last inequality uses the fact that, up to moving the cutting points
by at most $\lambda$, we can assume that the lengths involved are
multiples of $\lambda$, hence the factor $\ell/\lambda$ in front of the
max and the increase of the lengths by $2\lambda$.
\end{dem}

%\begin{prop}\label{cglocal}
%Suppose that, for some $A>1$, for any $A\lambda/4\leq \ell\leq A\lambda$ one has
%\[
%\abs{W_\ell}\leq (2m-1)^{\eta \ell}
%\]
%Then for any $\ell\geq A\lambda/4$,
%\[
%\abs{W_\ell}\leq (2m-1)^{(\eta+o(1)_{A\rightarrow\infty})\ell}
%\]
%where the constant implied in $o(1)$ depends only on $C$.
%\end{prop}
%
%\begin{dem}

Now for the proof of Proposition~\ref{cglocal} proper.

First, choose $\ell$ between $A\lambda$ and $4A\lambda/3$. By
Corollary~\ref{quasimult} and the assumptions, we have
\[
\abs{W_{\ell}}\leq (2m-1)^{\eta(\ell+4\alpha\lambda\log(\ell/C\lambda)
+6\lambda)+\log_{2m-1}(\ell/\lambda)}
\]

Let $B$ be a number (depending on $C$) such that
\[
4\alpha \log(B/C)+6+\frac1\eta \log_{2m-1}B \leq B
\]
(noting that $m\geq 2$, $\eta\geq 1/2$ and $\alpha\leq 1/C$
one can check that $B\geq 144/C^2$ is enough). It is then easy to check
that for $B'\geq B$ one has
\[
4\alpha \log(B'/C)+6+\frac1\eta \log_{2m-1}B' \leq 2\sqrt{B'B}
\]

Thus, if $\ell\geq A\lambda$ and $A\geq B$ we have
\[
\abs{W_{\ell}}\leq
(2m-1)^{\eta\left(\ell+2\lambda\sqrt{AB}\right)}\leq(2m-1)^{\eta\ell\left(1+2\sqrt{B/A}\right)}
\]

We have just shown that if $\abs{W_\ell}\leq (2m-1)^{\eta\ell}$ for
$\ell\leq A\lambda$, then $\abs{W_\ell}\leq
(2m-1)^{\eta\ell\left(1+2\sqrt{B/A}\right)}$ for
$\ell\leq (4A/3)\lambda$. Thus, iterating the process shows that for
$\ell\leq (4/3)^kA\lambda$ we have \[
\abs{W_\ell}\leq (2m-1)^{\eta\ell\prod_{0\leq i<
k}\left(1+2\sqrt{\frac{B}{A}}\left(\frac34\right)^{i/2}\right)}
\]
and we are done as the product
$\prod_i\left(1+2\sqrt{\frac{B}{A}}\left(\frac34\right)^{i/2}\right)$
converges to some value tending to $1$ when $A\rightarrow\infty$; if one cares, its value is less than $\exp
\frac{200}{C\sqrt{A}}$.
\end{dem}

\section{Application to random groups: the free case}

Here we first treat the case when the initial group $G$ is the free group
$F_m$ on $m$ generators. This will serve as a template for the more
complex general case.

So let $G=\presgroup{a_1,\ldots,a_m}{R}$ be a random group at density
$d$, with $R$ a set of $(2m-1)^{d\ell}$ random reduced words.

We have to evaluate the number of reduced words of a given length $L$ which
represent the trivial element in $G$. Any such word is the boundary word
of some van Kampen diagram $D$ with respect to the set of relators $R$.
We will proceed as follows: for any diagram $D$ involving $n$ relators,
we will evaluate the expected number of $n$-tuples of random relators
from $R$ that make it a van Kampen diagram. We will show that this
expected number is controlled by the boundary length $L$ of the diagram,
and this will finally allow to control the number of van Kampen diagrams
of boundary length $L$.

We call a van Kampen diagram \emph{non-filamenteous} if each of its
edges lies on the boundary on some face. Each diagram can be decomposed
into non-filamenteous components linked by filaments. For the
filamenteous part we will use the estimation from~\cite{Ch93}, one step of
which counts the number of ways in which the different non-filamenteous
parts can be glued together to form a van Kampen diagram.

So we first focus on non-filamenteous diagrams, for which a genuinely new
argument has to be given compared to~\cite{Ch93} (since the number of
relators here is unbounded).

We first suppose that we care only about diagrams with at most $K$ faces,
for some $K$ to be chosen later. (We will of course use the locality of
cogrowth principle to remove this assumption.)

\subsection{Fulfilling of diagrams}

So let $D$ be a non-filamenteous van Kampen diagram. Let $\abs{D}$ be its
number of faces and let $n\leq K$ be the number of different relators it
involves. Let $m_i$, $1\leq i \leq n$ be the number of times the $i$-th
relator appears in $D$, where we choose to enumerate the relators in
decreasing order of multiplicity, that is, $m_1\geq m_2 \geq \ldots \geq
m_n$. Let also $D_i$ be the subdiagram of $D$ made of relators
$1,2,\ldots,i$ only, so that $D=D_n$.

It is shown in~\cite{Oll} (section 2.2) that to this
diagram we can associate numbers $d_1,\ldots,d_n$ such that
\begin{itemize}
\item The probability that $i$ given random relators fulfill
$D_i$ is less than $(2m-1)^{d_i-id\ell}$ ; consequently, the probability
that there exists an $i$-tuple of relators in $R$ fulfilling $D_i$ is
less than $(2m-1)^{d_i}$.
\item The following isoperimetric inequality holds :
\[
\abs{\d D}\geq (1-2d)\ell\abs{D}+2\sum d_i(m_i-m_{i+1})
\]
\end{itemize}

So for any fixed $\eps>0$ we can suppose that $d_i\geq -\eps\ell$ for all $i$
(otherwise, $D$ appears as a van Kampen diagram of the random
presentation with probability less than $(2m-1)^{-\eps \ell}$, which
tends exponentially to $0$ as $\ell\rightarrow\infty$).

\def\myonehalf{{\textstyle \frac12}}

Then, using that $m_i-m_{i+1}\geq 0$ we can write
\begin{eqnarray*}
\abs{\d D} &\geq &
(1-2d)\ell\abs{D}-2\eps\ell\sum (m_i-m_{i+1})+2\sum
(d_i+\eps\ell)(m_i-m_{i+1})
\\&\geq &
(1-2d)\ell\abs{D}-2\eps\ell m_1 + 2 (d_n+\eps\ell)(m_n-m_{n+1})
\\&\geq &
\myonehalf (1-2d)\ell\abs{D} + 2 d_n
\end{eqnarray*}
where we chose to set $\eps=(1-2d)/4$ and where we used $m_1\leq
\abs{D}$ and $m_n\geq1$, $m_{n+1}=0$ by definition.

Now we know that for a given $n$-tuple of random relators, the
probability that this $n$-tuple fulfills $D$ is at most
$(2m-1)^{d_n-nd\ell}$. So, as there are $(2m-1)^{nd\ell}$ $n$-tuples of
relators in $R$, the expected number $S$ of $n$-tuples fulfilling $D$ in
$R$ is at most $(2m-1)^{d_n}$, which so turns out to be not only an upper
bound for the probability of $D$ to be fulfillable but rather an estimate
of the number of ways in which it is.  (The probabilities that two
$n$-tuples fulfill the diagram are independent only when the $n$-tuples
are disjoint, but expectation is linear anyway.)

By Markov's inequality, the probability to pick a random presentation $R$
for which $S\geq (2m-1)^{\eps''\ell}(2m-1)^{d_n}$ is less than
$(2m-1)^{-\eps''\ell}$.

Thus, for any fixed integer $K$ and any $\eps''>0$, with probability
exponentially close to $1$ as $\ell\rightarrow\infty$, we can suppose
that a given (hence any, since the number of diagrams with less than $K$
faces grows subexponentially) non-filamenteous diagram can be filled in
at most $(2m-1)^{\eps''\ell}(2m-1)^{d_n}$ different ways by relators of
$R$. (The $\ell$ up from which this holds depends of course on $\eps''$
and $K$.)

The last inequality above can be rewritten as
\[
d_n\leq \frac{1}{2}\left(\abs{\d D}-(\myonehalf-d)\ell\abs{D}\right)
\]
or as $\abs{D}\geq 1$
\[
d_n+\eps''\ell\leq \frac{1}{2}\left(\abs{\d D}-(\myonehalf-d-2\eps'')\ell\abs{D}\right)
\]
so if we choose $\eps''\leq(\myonehalf-d)/2$, this is at most $\abs{\d D}/2$.

The conclusion is:

\begin{prop}
\label{numbernonfil}
For each $K$, with
probability exponentially close to $1$ as $\ell\rightarrow\infty$,
for each
non-filamenteous van Kampen diagram with at most $K$ faces, the number of ways to fulfill it with
relators of $R$ is at most $(2m-1)^{\abs{\d D}/2}$.
\end{prop}

\subsection{Evaluation of the cogrowth}

We now conclude using the general scheme of~\cite{Ch93}, together
with Proposition~\ref{cglocal} which enables to check only a finite
number of diagrams.

Consider a reduced word $w$ in the generators $a_i^{\pm
1}$, representing $e$ in the random group. This word is the boundary word
of some van Kampen diagram $D$ which may have filaments.

Choose $\eps>0$. We are going to show that with probability exponentially
close to $1$ when $\ell\rightarrow \infty$, the number of such words $w$
is at most $(2m-1)^{(1/2+\eps)\abs{w}}$.

\bigskip

We know from~\cite{Oll} (Section~2.2) that up to exponentially small
probability in $\ell$, we
can suppose that any
diagram satisfies the inequality
\[
\abs{\d D}\geq C\ell\abs{D}
\]
where $C$ depends only on the density $d$ (basically $C=1/2-d$ divided by
the constants appearing in the Cartan-Hadamard-Gromov theorem,
see~\cite{Oll}) and not on $\ell$.

Now we use Proposition~\ref{cglocal}. We are facing a group $G$ in which
all relations are of length $\ell$. Consider a constant $A$ given by
Proposition~\ref{cglocal} such that if we know that $\abs{W_L} \leq
(2m-1)^{L(1/2+\eps/2)}$ for $L\leq A\ell$, then we know that $\abs{W_L}\leq
(2m-1)^{L(1/2+\eps)}$ for any $L$. Such an $A$ depends only on the
isoperimetry constant $C$.

So we suppose that our word $w$ has length at most $A\ell$. We have
$\abs{w}=\abs{\d D}\geq C\ell\abs{D}$ and in particular, $\abs{D}\leq
A/C$, which is to say, we have to consider only diagrams with a number of
faces bounded independently of $\ell$.

So set $K=A/C$, which most importantly does not depend on $\ell$. After
Proposition~\ref{numbernonfil}, we can assume (up to exponentially small
probability) that for any non-filamenteous diagram $D'$ with at most $K$
faces, the number of ways to fulfill it with relators of the random
presentation is at most $(2m-1)^{\abs{\d D'}/2}$.

\bigskip

Back to our word $w$ read on the boundary of some diagram $D$. Decompose
$D$ into filaments and connected non-filamenteous parts $D_i$. The word
$w$ is determined by the following data: a set of relators from the
random presentation $R$ fulfilling the $D_i$'s, a set of reduced
words to put on the filaments, the combinatorial choice of the diagrams
$D_i$, and the combinatorial choice of how to connect the $D_i$'s using
the filaments.

The combinatorial part is precisely the one analyzed in~\cite{Ch93}. It
is shown there (section ``Premier pas'') that if each $D_i$ satisfies
$\abs{\d D_i}\geq L$, the combinatorial factor controlling the connecting
of the $D_i$'s by the filaments and the sharing of the length $\abs{\d
D}$ between the filaments and the $D_i$'s is less than
\[
\frac{\abs{w}}{L} \abs{w} (eL)^{2\abs{w}/L} (2eL)^{\abs{w}/L}
(3eL)^{2\abs{w}/L}
\]
Observe that for $L$ large enough this behaves like
$(2m-1)^{\abs{w}O(\log L/L)}$.

Here each diagram $D_i$ satisfies $\abs{\d D_i}\geq C\ell\abs{D_i}\geq
C\ell$, so setting $L=C\ell$, each $D_i$ has boundary length at
least $L$. In particular, $O(\log L/L)=O(\log \ell/\ell)$.

The number of components $D_i$ is obviously at most $\abs{w}/L$. Each
component has at most $K$ faces since $D$ itself has. So the number of
choices for the combinatorial choices of the diagrams $D_i$'s is at most
$N(K)^{\abs{w}/L}$ where $N(K)$ is the (finite!) number of planar graphs
with at most $K$ faces. This behaves like $(2m-1)^{\abs{w}O(1/L)}$.

Now the number of ways to fill the $D_i$'s with relators from the random
presentation is, after Proposition~\ref{numbernonfil}, at most $\prod
(2m-1)^{\abs{\d D_i}/2}= (2m-1)^{\sum\abs{\d D_i}/2}$.

The last choice to take into account is the choice of reduced words to
put on the filaments.  The total length of the filaments is
$\frac12 (\abs{w}-\sum\abs{\d D_i})$ (each edge of a filament counts twice in
the boundary), thus the number of ways to fill in the filaments is at
most $(2m-1)^{\frac12\left(\abs{w}-\sum\abs{\d D_i}\right)}$.

So the total number of possibilities for $w$ is
\[
(2m-1)^{\abs{w}O(\log\ell/\ell)+\frac12\left(\abs{w}-\sum\abs{\d
D_i}\right)+\sum\abs{\d D_i}/2}
\]
and if we take $\ell$ large enough, this will be at most
$(2m-1)^{\abs{w}(1/2+\eps/2)}$, after what we conclude by
Proposition~\ref{cglocal}.

This proves Theorem~\ref{cormain}.

\section{The non-free case}
\label{nonfree}

Now we deal with random quotients of a non-elementary torsion-free
hyperbolic group $G_0$.  We are going to give the proof in the case of a
random quotient by plain random words, the case of a quotient by random
reduced words being similar.

So let $G_0$ be a non-elementary torsion-free hyperbolic group given by
the presentation $\presgroup{a_1,\ldots,a_m}{Q}$ ($m\geq 2$), with the
relations in $Q$ having length at most $\lambda$. Let $\theta$ be the
gross cogrowth of $G_0$ w.r.t.\ this generating set. Let $G=G_0/\langle
R\rangle$ be a random quotient of $G_0$ by a set $R$ of $(2m)^{d\ell}$
randomly chosen words of length $\ell$.  Also set $\beta=1-\theta$, so
that the random quotient axioms of~\cite{Oll} (section~4) are satisfied.

\bigskip

We have to show that the number of boundary words of van Kampen diagrams
of a given boundary length $L$ grows slower than $(2m)^{L(\theta+\eps)}$.
This time, since we are going to give a proof in the case of gross
cogrowth rather than cogrowth, we will not have many problems with
filaments: the counting of filaments is already included in the knowledge of
gross cogrowth of $G_0$.

For a van Kampen diagram $D$, let $D''$ be the subdiagram made of faces
bearing ``new'' relators in $R$, and $D'$ be the part made of faces
bearing ``old'' relators in $Q$. By Proposition~32
of~\cite{Oll}, we know that very probably $G$ is hyperbolic and that its
isoperimetric inequality takes the form
\[
\abs{\d D}\geq \kappa\ell\abs{D''}+\kappa'\abs{D'}
\]
whenever $D$ is reduced and $D'$ is minimal, with $\kappa,\kappa'>0$ and
where, most importantly, $\kappa$ and $\kappa'$ do not depend on $\ell$.
By definition of $\A(D)$, this can be rewritten as $\abs{\d D}\geq C \A(D)$ with $C=\min(\kappa,
\kappa'/\lambda)$.

Fix some $\eps>0$ and let $A$ be the constant provided by
Proposition~\ref{cglocal} applied to $G$, having the property that if we
know that gross cogrowth is at most $\theta+\eps/2$ up to words of length
$A\ell$, then we know that gross cogrowth is at most $\theta+\eps$. This
$A$ depends on $\eps$, $C$ and $G_0$ but not on $\ell$. Thanks to this and the
isoperimetric inequality, we only have to consider diagrams of boundary
length at most $A\ell$ hence area at
most $A\ell/C$. In particular the number of new relators $\abs{D''}$ is
at most $A/C$. So for all the sequel set
\[
K=A/C
\]
which, most importantly, does not depend on $\ell$. This is the maximal
size of diagrams we have to consider, thanks to the local-global
principle.

\subsection{Reminder from~\cite{Oll}}

In this context, it is proven in~\cite{Oll} that the van Kampen diagram
$D$ can be seen as a ``van Kampen diagram at scale $\ell$ with respect to
the new relators, with equalities modulo $G_0$''. More precisely, this
can be stated as follows: (we refer to~\cite{Oll} for the definition of
``strongly reduced'' diagrams; the only thing to know here is that for
any word equal to $e$ in $G$, there exists a strongly reduced van Kampen
diagram with this word as its boundary word).

\begin{prop}[ (\cite{Oll}, section~6.6)]
\label{propdavKd}
Let $G_0=\presgroup{S}{Q}$ be a non-elementary hyperbolic group, let $R$ be a set of
words of length $\ell$, and consider the group $G=G_0/\langle R
\rangle=\presgroup{S}{Q \cup R}$.

Let $K\geq 1$ be an arbitrarily large integer and let $\eps_1,\eps_2>0$
be arbitrarily small numbers. Take $\ell$ large enough depending on $G_0,
K,\eps_1,\eps_2$.

Let $D$ be a van Kampen diagram with respect to the presentation
$\presgroup{S}{Q \cup R}$, which is strongly reduced, of area at most
$K\ell$. Let also $D'$ be the subdiagram of $D$ which is the union of the
$1$-skeleton of $D$ and of those faces of $D$ bearing relators in $Q$
(so $D'$ is a possibly non-simply connected van Kampen
diagram with respect to $G_0$), and suppose that
$D'$ is minimal.

We will call \emph{worth-considering} such a van Kampen diagram.

Let $w_1,\ldots,w_p$ be the boundary (cyclic) words of $D'$, so that each
$w_i$ is either the boundary word of $D$ or a relator in $R$.

\smallskip

Then there exists an integer $k\leq 3K/\eps_2$ and words
$x_2,\ldots,x_{2k+1}$ such that:
\begin{itemize}
\item Each $x_i$ is a subword of some cyclic word $w_j$;

\item As subwords of the $w_j$'s, the $x_i$'s are disjoint and their
union exhausts a proportion at least $1-\eps_1$ of the total length of
the $w_j$'s.

\item For each $i\leq k$, there exists words $\delta_1,\delta_2$ of
length at most $\eps_2 (\abs{x_{2i}}+\abs{x_{2i+1}})$ such that
$x_{2i}\delta_1 x_{2i+1} \delta_2=e$ in $G_0$.

\item If two words $x_{2i}$, $x_{2i+1}$ are subwords of the boundary words of
two faces of $D$ bearing the same relator $r^{\pm 1}\in R$, then, as
subwords of $r$, $x_{2i}$ and $x_{2i+1}$ are either disjoint or equal with
opposite orientations (so that the above equality reads $x\delta_1
x^{-1}\delta_2=e$).
\end{itemize}

The couples $(x_{2i},x_{2i+1})$ are called \emph{translators}.
Translators are called \emph{internal}, \emph{internal-boundary} or
\emph{boundary-boundary} according to whether $x_{2i}$ and $x_{2i+1}$
is a subword of some $w_j$ which is a relator in $R$ or the boundary word of
$D$.
\end{prop}

(There are slight differences between the presentation here and that
in~\cite{Oll}. Therein, boundary-boundary translators did not have to be
considered: they were eliminated earlier in the process, before
section~6.6, because they have a positive contribution to boundary
length, hence always improve isoperimetry and do not deserve
consideration in order to prove hyperbolicity. Moreover, in~\cite{Oll} we
further distinguished ``commutation translators'' for the kind of
internal translator with $x_{2i}=x_{2i+1}^{-1}$, which we need not do
here.)

Translators appear as dark strips on the following figure:

\begin{center}
\includegraphics{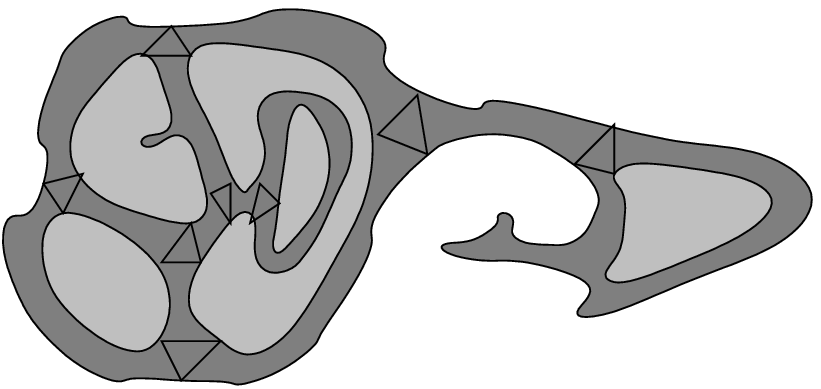}
\end{center}

\begin{rem}
\label{numberofdavKds}
The number of ways to partition the words
$w_i$ into translators is at most $(2K\ell)^{12K/\eps_2}$, because each
$w_i$ can be determined by its starting- and endpoint, which can be given
as numbers between $1$ and $2K\ell$ which is an upper bound for the
cumulated length of the $w_i$'s (since the area of $D$ is at most
$K\ell$). For fixed $K$ and $\eps_2$ this grows subexponentially in
$\ell$.
\end{rem}

\begin{rem}
\label{orderboundary}
Knowing the words $x_i$, the number of possibilities for the boundary
word of the diagram is at most $(6K/\eps_2)!$ (choose which subwords
$x_i$ make the boundary word of the diagram, in which order), which does
not depend on $\ell$ for fixed $K$ and $\eps_2$.
\end{rem}

\bigskip

We need another notion from~\cite{Oll}, namely, that of \emph{apparent
length} of an element in $G_0$. This basically answers the question: If
this element were obtained through a random walk at time $t$, what would
be a reasonable value of $t$? This accounts for the fact that, unlike in
the free group, the hitting probability of an element in the group does
not depend only on the norm of this element.

Apparent length is defined in~\cite{Oll} in a more general setting, with
respect to a measure on the group, which is here the measure obtained
after a simple random walk with respect to the given set of generators
$a_1,\ldots,a_m$. We only give here what the definition amounts to in our
context.

\def\L{\mathbb{L}}

\begin{defi}[ (Definition~36 of~\cite{Oll})]
Let $x$ be a word. Let $\eps_2>0$. Let $L$ be an integer. Let $p_L(xuyv=e)$
be the probability that, for a random word $y$ of length $L$, there
exists elements $u,v\in G_0$ of norm at most $\eps_2(\abs{x}+L)$ such that
$xuyv=e$ in $G_0$.

The \emph{apparent length} of $x$ at test-length $L$ is
\[
\L_L(x)=-\frac{1}{1-\theta} \log_{2m} p_L(xuyv=e) - L
\]

The \emph{apparent length} of $x$ is
\[
\L(x)=\min\left(\norm{x}\frac{\theta}{1-\theta}\,,\;\min_{0 \leq L
\leq K\ell} \L_L(x)\right)
\]
where we recall $\ell$ is the length of the relators in a random
presentation.
\end{defi}

(The first term $\norm{x}\theta/(1-\theta)$ is an easy upper bound for
$\ell_{\norm{x}}(x)$, and so if $\norm{x}\leq K\ell$ then the first term
in the min is useless.)

It is shown in~\cite{Oll}, section~6.7, that in a randomly chosen
presentation at density $d$ and length $\ell$, all subwords of the
relators have apparent length at most $4\ell$, with probability
exponentially close to $1$ as $\ell\rightarrow\infty$. So from now on we
suppose that this is indeed the case.

%We are going to consider apparent lengths of subwords of the random
%relators in the presentation. We need that these apparent lengths be not
%too small, namely, of the same order of magnitude as the usual length.
%Precisely:
%
%\begin{prop}[ (Proposition~(REFERENCE) of~\cite{Oll})]
%Let $x$ be a random word of length $\abs{x}$. Then,
%uniformly for any $\L\leq \abs{x}$ we have
%\[
%\Pr\left(\L(X)\leq \abs{x}-L\right)\leq
%(2m)^{-(1-\theta)L+o(\abs{x})_{\abs{x}\rightarrow\infty}}
%\]
%where the $o(\abs{x})$ depends only on $G_0$.
%\end{prop}

\bigskip

We further need the notion of a \emph{decorated abstract van Kampen
diagram} (which was implicitly present in the free case when we mentioned
the probability that some diagram ``is fulfilled by random relators''),
which is inspired by Proposition~\ref{propdavKd}: it carries the
combinatorial information about how the relators and boundary word of a
diagram were cut into subwords in order to make the translators.

\def\D{\mathcal{D}}

\begin{defi}[ (Decorated abstract van Kampen diagram)]
\label{defdavKds}
Let $K\geq 1$ be an arbitrarily large integer and let $\eps_1,\eps_2>0$
be arbitrarily small numbers. Let $I_\ell$ be the cyclically ordered set of
$\ell$ elements.

A \emph{decorated abstract van Kampen diagram} $\D$ is the following data:
\begin{itemize}
\item An integer $\abs{\D}\leq K$ called its \emph{number of faces}.
\item An integer $\abs{\d \D}\leq K\ell$ called its \emph{boundary length}.
\item An integer $n\leq \abs{\D}$ called its \emph{number of distinct
relators}.
\item An application $r^\D$ from $\{1,\ldots,\abs{\D}\}$ to
$\{1,\ldots,n\}$; if $r^\D(i)=r^\D(j)$ we will say that \emph{faces $i$
and $j$ bear the same relator}.
\item An integer $k\leq 3K/\eps_2$ called the \emph{number of
translators} of $\D$.
\item For each integer $2\leq i\leq 2k+1$, a set of the form $\{j_i\}\times
I'_i$ where either $j_i$ is an integer between $1$ and $\abs{\D}$ and
$I'_i$ is an oriented cyclic subinterval of
$I_\ell$, or $j_i=\abs{\D}+1$ and $I'_i$ is a subinterval of $I_{\abs{\d
\D}}$;
this is called an
\emph{(internal) subword of the $j_i$-th face} in the first case, or a
\emph{boundary subword} in the second case.
\item For each integer $1\leq i \leq k$ such that $j_{2i}\leq \abs{\D}$, an integer between $0$ and
$4\ell$ called the \emph{apparent length of the $2i$-th subword}.
\end{itemize}
such that
\begin{itemize}
\item The sets $\{j_i\}\times I'_i$ are all disjoint and the cardinal
of their union is at least $(1-\eps_1)\left(\abs{\D}\ell+\abs{\d \D}\right)$.
\item For all $1\leq i\leq k$ we have $j_{2i}\leq j_{2i+1}$ (this can be
ensured by maybe swapping them).
\item If two faces $j_{2i}$ and $j_{2i+1}$ bear the same relator, then
either $I'_{2i}$ and $I'_{2i+1}$ are disjoint or are equal with opposite
orientations.
\end{itemize}
\end{defi}

This way, Proposition~\ref{propdavKd} ensures that any
worth-considering van Kampen diagram $D$ with respect to $G_0/\langle R
\rangle$ defines a decorated abstract van Kampen diagram $\D$ in the way
suggested by terminology (up to rounding the apparent lengths to the
nearest integer; we neglect this problem). We will say that $\D$ is
\emph{associated to} $D$.  Remark~\ref{numberofdavKds} tells that the
number of decorated abstract van Kampen diagrams grows subexponentially
with $\ell$ (for fixed $K$).

Given a decorated abstract van Kampen diagram $\D$ and $n$ given relators
$r_1,\ldots,r_n$, we say that these relators \emph{fulfill $\D$} if there
exists a worth-considering van Kampen diagram $D$ with respect to
$G_0/\langle r_1,\ldots,r_n\rangle$, such that the associated decorated
abstract van Kampen diagram is $\D$. Intuitively speaking, the relators
$r_1,\ldots,r_n$ can be ``glued modulo $G_0$ in the way described by
$\D$''.

So we want to study which diagrams can probably be fulfilled by random
relators in $R$. The main conclusion from~\cite{Oll} is that these are
those with large boundary length, hence hyperbolicity of the quotient
$G_0/\langle R\rangle$. Here for cogrowth we are rather interested in
the number of ways to fulfill an abstract diagram with given boundary
length.

\subsection{Cogrowth of random quotients}

So now let $R$ again be a set of $(2m)^{d\ell}$ random relators. Let $\D$
be a given decorated abstract van Kampen diagram. Recall we set $K=A/C$.
The free parameters $\eps_1$ and $\eps_2$ will be chosen later.

We will show (Proposition~\ref{lastprop}) that, up to exponentially small
probability in $\ell$, the number of different boundary words of
worth-considering van Kampen diagrams $D$ such that $\D$ is associated
to $D$, is at most $(2m)^{\theta\abs{\d \D}(1+\eps/2)}$.

\paragraph{Further notations.}
Let $n$ be the number of distinct relators in $\D$. For $1\leq a\leq n$,
let $m_a$ be the number of times the $a$-th relator appears in $\D$. Up
to reordering, we can suppose that the $m_a$ 's are non-increasing. Also
to avoid trivialities take $n$ minimal so that $m_n\geq 1$.

Let also $P_a$ be the probability that, if $a$ words $r_1,\ldots,r_a$ of
length $\ell$ are picked at random, there exist $n-a$ words
$r_{a+1},\ldots,r_n$ of lengt $\ell$ such that the relators
$r_1,\ldots,r_n$ fulfill $\D$. The $P_a$ 's are of course a
non-increasing sequence of probabilities. In particular, $P_n$ is the
probability that a random $n$-tuple of relators fulfills $\D$.

Back to our set $R$ of $(2m)^{d\ell}$ randomly chosen relators. Let $P^a$
be the probability that there exist $a$ relators $r_1,\ldots,r_a$ in
$R$, such that there exist words $r_{a+1}, \ldots,r_n$ of length $\ell$
such that the relators $r_1,\ldots,r_n$ fulfill $\D$. Again the $P^a$ 's
are a non-increasing sequence of probabilities and of course we have
\[
P^a\leq (2m)^{ad\ell} P_a
\]
since the $(2m)^{ad\ell}$ factor accounts for the choice of the $a$-tuple
of relators in $R$.

The probability that there exists a van Kampen diagram $D$ with respect
to the random presentation $R$, such that $\D$ is associated to $D$, is
by definition less than $P^a$ for any $a$. In particular, if for some
$\D$ we have $P^a\leq (2m)^{-\eps'\ell}$, then with probability
exponentially close to $1$ when $\ell\rightarrow\infty$, $\D$ is not
associated to any van Kampen diagram of the random presentation. Since,
by Remark~\ref{numberofdavKds}, the number of possibilities for $\D$
grows subexponentially with $\ell$, we can sum this over $\D$ and
conclude that for any $\eps'>0$, with probability
exponentially close to $1$ when $\ell\rightarrow\infty$ (depending on
$\eps'$), all decorated abstract van Kampen diagrams $\D$ associated to
some van Kampen diagram of the random presentation satisfy $P^a\geq
(2m)^{-\eps'\ell}$ and in particular
\[
P_a \geq (2m)^{-ad\ell-\eps'\ell}
\]
which we assume from now on.

\medskip

We need to define one further quantity. Keep the notations of
Definition~\ref{defdavKds}. Let $1\leq a \leq n$ and let $1\leq i \leq k$
where $k$ is the number of translators of $\D$. Say that the $i$-th
translator is half finished at time $a$ if $r^\D(j_{2i})\leq a$ and
$r^\D(j_{2i+1})>a$, that is, if one side of the translator is a subword
of a relator $r_{a'}$ with $a'\leq a$ and the other of $r_{a''}$ with
$a''>a$. Now let $A_a$ be the sum of the apparent lengths of all
translators which are half finished at time $a$. In particular, $A_n$ is
the sum of the apparent lengths of all subwords $2i$ such that $2i$ is an
internal subword and $2i+1$ is a boundary subword of $\D$.

\paragraph{The proof.}
In this context, equation $(\star)$ (section~6.8) of~\cite{Oll} reads
\[
A_a-A_{a-1}\geq
m_a\left(\ell(1-\eps'')+\frac{\log_{2m}P_a-\log_{2m}P_{a-1}}{\beta}\right)
\]
where $\eps''$ tends to $0$ when our free parameters $\eps_1,\eps_2$ tend
to $0$ (and $\eps''$ also absorbs the $o(\ell)$ term in~\cite{Oll}). Also
recall that in the model of random quotient by plain random words, we
have
\[
\beta=1-\theta
\]
by Proposition~15 of~\cite{Oll}.

Setting $d'_a=\log_{2m}P_a$ and summing over $a$ we get, using $\sum
m_a=\abs{\D}$, that
\begin{eqnarray*}
A_n&\geq&
\left(\sum m_a\right)\ell \left(1-\eps''\right)+\frac1\beta \sum m_a(d'_a-d'_{a-1})
\\&=&\abs{\D}\ell(1-\eps'') +\frac1\beta \sum d'_a(m_a-m_{a+1})
\end{eqnarray*}

Now recall we saw above that for any $\eps'>0$, taking $\ell$ large enough
we can suppose that $P_a \geq (2m)^{-ad\ell-\eps'\ell}$, that is,
$d'_a+ad\ell+\eps'\ell\geq 0$. Hence
\begin{eqnarray*}
A_n&\geq&
\abs{\D}\ell(1-\eps'') + \frac1\beta \sum
(d'_a+ad\ell+\eps'\ell)(m_a-m_{a+1})
\\& &-\frac1\beta \sum
(ad\ell+\eps'\ell)(m_a-m_{a+1})
\\&=&
\abs{\D}\ell(1-\eps'') +\frac1\beta
\sum(d'_a+ad\ell+\eps'\ell)(m_a-m_{a+1})-\frac{d\ell}{\beta}\sum
m_a-\frac{\eps'\ell}{\beta} m_1
\\&\geq& \abs{\D}\ell(1-\eps'') +\frac{d'_n+nd\ell+\eps'\ell}{\beta}m_n
-\frac{d\ell+\eps'\ell}{\beta}\sum m_a
\end{eqnarray*}
where the last inequality follows from the fact that we chose the order
of the relators so that $m_a-m_{a+1}\geq 0$.

So using $m_n\geq 1$ we finally get
\[
A_n\geq
\abs{\D}\ell\left(1-\eps''-\frac{d+\eps'}{\beta}\right)+\frac{d'_n+nd\ell}{\beta}
\]

Suppose the free parameters $\eps_1$, $\eps_2$ and $\eps'$ are chosen
small enough so that $1-\eps''-(d+\eps)/\beta\geq 0$ (remember that
$\eps''$ is a function of $\eps_1,\eps_2$ and $K$; we will further
decrease $\eps_1$ and $\eps_2$ later). This is possible
since by assumption we take the density $d$ to be less than the critical
density $\beta$. This is the only, but crucial, place where density plays
a role. Thus the first term in the inequality above is non-negative and
we obtain the simple inequality $A_n\geq (d'_n+nd\ell)/\beta$.

\begin{prop}
\label{alntuplesineq}
Up to exponentially small probability in $\ell$, we can suppose that any
decorated abstract van Kampen diagram $\D$ satisfies
\[
A_n(\D)\geq\frac{d'_n(\D)+nd\ell}{\beta}
\]
\end{prop}

\bigskip

This we now use to evaluate the number of possible boundary words for
van Kampen diagrams associated with $\abs{\D}$.

Remember that, by definition, $d'_n$ is the log-probability that $n$ random
relators $r_1,\ldots,r_n$ fulfill $\D$. As there are $(2m)^{nd\ell}$
$n$-tuples of random relators in $R$ (by definition of the density
model), by linearity of expectation the expected number of $n$-tuples of
relators in $R$ fulfilling $\D$ is $(2m)^{nd\ell+d'_n}$, hence the
interest of an upper bound for $d'_n+nd\ell$.

By the Markov inequality, for given $\D$ the probability to pick a random
set $R$ such that the number of $n$-tuples of relators of $R$ fulfilling
$\D$ is greater than $(2m)^{nd\ell+d'_n+C\eps\ell/4}$, is less than
$(2m)^{-C\eps\ell/4}$. By Remark~\ref{numberofdavKds} the number of
possibilities for $\D$ is subexponential in $\ell$, and so, using
Proposition~\ref{alntuplesineq} we get

\begin{prop}
\label{numberntuples}
Up to exponentially small probability in $\ell$, we can suppose that for
any decorated abstract van Kampen diagram $\D$, the number of $n$-tuples
of relators in $R$ fulfilling $\D$ is at most 
\[
(2m)^{\beta A_n(\D)+ C\eps\ell/4}
\]
\end{prop}

Now let $D$ be a van Kampen diagram associated to $\D$. Given $\D$ we want to
evaluate the number of different boundary words for $\D$. Recall
Proposition~\ref{propdavKd}: the boundary word of $D$ is determined by
giving two words for each boundary-boundary translator, and one word for
each internal-boundary translator, this last one being subject to the
apparent length condition imposed in the definition of $\D$. By
Remark~\ref{orderboundary}, the number of ways to combine these subwords
into a boundary word for $D$ is controlled by $K$ and $\eps_2$
(independently of $\ell$).

So let $(x_{2i},x_{2i+1})$ be a boundary-boundary translator in $D$. By
Proposition~\ref{propdavKd} (definition of translators) there exist words
$\delta_1,\delta_2$ of length at most
$\eps_2(\abs{x_{2i}}+\abs{x_{2i+1}})$ such that
$x_{2i}\delta_1x_{2i+1}\delta_2=e$ in $G_0$. So
$x_{2i}\delta_1x_{2i+1}\delta_2$ is a word representing the trivial
element in $G_0$, and by definition of $\theta$ the number of
possibilities for $(x_{2i},x_{2i+1})$ is at most
$(2m)^{\theta(\abs{x_{2i}}+\abs{x_{2i+1}})(1+2\eps_2)}$.

Now let $(x_{2i},x_{2i+1})$ be an internal-boundary translator.
The apparent length of $x_{2i}$ is imposed in the definition of $\D$. The
subword $x_{2i}$ is an internal subword of $D$, and so by definition is a
subword of some relator $r_i\in R$. So if the relators in $D$ are given,
$x_{2i}$ is determined. But knowing $x_{2i}$ still leaves open lots of
possibilities for $x_{2i+1}$. This is where apparent length comes into
play.

Since $y=x_{2i+1}$ is a boundary word of $D$ one has $\abs{y}\leq
A\ell\leq K\ell$. So by definition we have $\L(x)\leq
\L_{\abs{y}}(x_{2i})$. By definition of translators there exist words $u$
and $v$ of length at most $\eps_2\ell$ such that $x_{2i}uyv=e$ in $G_0$.
By definition of $\L_{\abs{y}}(x_{2i})$, if $y'$ is a random word of
length $\abs{y}$, then the probability that $x_{2i}uy'v=e$ in $G_0$ is
$(2m)^{-(1-\theta)\left(\abs{y}+\L_{\abs{y}}(x_{2i})\right)}\leq
(2m)^{-(1-\theta)\left(\abs{y}+\L(x_{2i})\right)}$. This means that the
total number of words $y'$ of length $\abs{y}$ such that there exists
$u$, $v$ with $x_{2i}uyv=e$ is at most
$(2m)^{\abs{y}}(2m)^{-(1-\theta)\left(\abs{y}+\L(x_{2i})\right)}=(2m)^{\theta\abs{y}-(1-\theta)\L(x_{2i})}$.
So, given $x_{2i}$, the number of possibilities for $y=x_{2i+1}$ is less
than this number.

So if the relators in $R$ fulfilling $\D$ are fixed, the number of
possible boundary words for $D$ is the product of
$(2m)^{\theta(\abs{x_{2i}}+\abs{x_{2i+1}})(1+2\eps_2)}$ for all
boundary-boundary translators $(x_{2i},x_{2i+1})$, times the product of
$(2m)^{\theta\abs{x_{2i+1}}-(1-\theta)\L(x_{2i})}$ for all internal-boundary
translators $(x_{2i},x_{2i+1})$, times the number of ways to order these
subwords (which is subexponential in $\ell$ by
Remark~\ref{orderboundary}), times the number of possibilities for the
parts of the boundary of $D$ not belonging to any translator, which by
Proposition~\ref{propdavKd} have total length not exceeding
$\eps_1K\ell$.

Now the sum of $\abs{x_{2i}}+\abs{x_{2i+1}}$ for all boundary-boundary
translators $(x_{2i},x_{2i+1})$, plus the sum of $\abs{x_{2i+1}}$ for all
internal-boundary translators, is $\abs{\d \D}$ (maybe up to
$\eps_1K\ell$). And the sum of $\L(x_{2i})$ for all internal-boundary
translators is $A_n$ by definition.

So given $\D$ and given a $n$-tuple of relators fulfilling $\D$, the
number of possibilities for the boundary word of $D$ is at most
\[
(2m)^{\theta\abs{\d \D}(1+2\eps_2)-(1-\theta)A_n+\eps_1K\ell}
\]
up to a subexponential term in $\ell$.
By Proposition~\ref{numberntuples}
(remember $\beta=1-\theta$),
if we include the choices of the relators fulfilling $\D$ the number of
possibilities is at most 
\[
(2m)^{\theta\abs{\d
\D}(1+2\eps_2)+\eps_1K\ell+C\eps\ell/4}
\]

If we choose $\eps_2\leq\eps/16$ and $\eps_1\leq\eps C/8K$ so that (using
$\abs{\d D}\geq C\ell\abs{\D}\geq C\ell$ for any fulfillable abstract
diagram) the sum of the corresponding terms is less than $\eps\abs{\d
D}/4$ (note that this choice does not depend on $\ell$) and if we
remember that, after Remark~\ref{numberofdavKds}, the number of choices
for $\D$ is subexponential in $\ell$, we finally get:

\begin{prop}
\label{lastprop}
Up to exponentially small probability in $\ell$, the number of different
boundary words of worth-considering van Kampen diagrams of a random
presentation with given boundary length $L$, is at most
\[
(2m)^{\theta L (1+\eps/2)}
\]
\end{prop}

But remember the discussion at the beginning of section~\ref{nonfree}
(where we invoked Proposition~\ref{cglocal}): it is enough to show that gross
cogrowth is at most $\theta+\eps/2$ for words of length $L$ between
$A\ell/4$ and $A\ell$.  Any such word is the boundary word of a van
Kampen diagram of area at most $K\ell$, hence is the boundary word of
some worth-considering van Kampen diagram. This ends the proof of
Theorem~\ref{thmmain}.

\end{document}